\RequirePackage{fix-cm}
\documentclass[onecolumn]{svjour3questa}       
\usepackage{amsmath,amssymb}
\usepackage{mathptmx}      
\usepackage{hyperref}

\begin{document}

\title{Conjectures on Symmetric Queues in Heavy Traffic}
\author{ Bert Zwart}
\institute{CWI Amsterdam and Eindhoven University of Technology
\email{bert.zwart@cwi.nl}}

\date{\today}

\maketitle

\section{Introduction}

We consider the queue-length process 
in the $M/G/1$ queue with symmetric service discipline, which is defined as follows: with $n$ customers in the system, the server
works on the customer in position $i$ at rate $\gamma(n,i)\geq 0$. We assume $\sum_i\gamma(n,i)=1$  for all $n\geq 1$, to make the service discipline work conserving.  
The customer arrival process is Poisson with rate $\lambda$. A customer arriving to a queue of size $n$ chooses position $i$ with probability $\gamma(n + 1,i)$, moving each customer in position $k\geq i$ to position $k+1$. 
Conversely, when a customer at position $i$ departs, each customer in position $k>i$ moves to position $k-1$.  We refer to this system as a symmetric queue governed by $\gamma$.
Important special cases of symmetric service disciplines are  
Last-Come-First-Served (LCFS), where $\gamma(n,1) = 1$ for all $n\geq 1$, and Processor Sharing (PS), which is modeled by taking $\gamma (n,i)=1/n$ for all $1\leq i \leq n$ and all $n \geq 1$. We refer to \cite{Kelly1979} for more background. 

Symmetric service disciplines have the appealing property that the stationary distribution of the queue-length process $(Q^\gamma(t), t \geq 0)$, if it exists, is {\em insensitive} to the service-time distribution, apart from its mean $m$. 
In particular if the load $\rho = \lambda m<1$, then, as $t\rightarrow\infty$, $Q^\gamma(t)$ converges weakly to a random variable $Q^\gamma(\infty)$ which satisfies 
\begin{equation}
    P(Q^\gamma(\infty) \geq k) = \rho^k, k\geq 0.
\end{equation}
This note states two open problems related to the queue-length {\em process} in heavy traffic.
\section{Problem Statement}

Assume first that the second moment of the service-time distribution is finite. Let $r$ be a scaling parameter and let $\lambda_r$ be a sequence of arrival rates such that $\lambda_r \rightarrow 1/m$
and $r(1-\rho_r) = r(1-\lambda_r m)$ converges to a real-valued constant. Let 
$Q_r^\gamma(t), t\geq 0$ be the corresponding queue-length process in the symmetric  queue governed by $\gamma$ with Poisson arrival rate $\lambda_r$. For each $r$, 
define the re-scaled queue-length process $\hat Q_r^\gamma (t) = Q_r^\gamma (r^2t)/r, t\geq 0$.\\ 

\noindent
{\bf Problem 1:} show that $\hat Q_r^\gamma (\cdot)$ converges in distribution to a reflected Brownian motion (RBM) as $r\rightarrow\infty$. \\

\noindent
Our second problem relates to the case where the variance of the service-time distribution is {\em infinite}. 
To this end, let $L$ be a slowly varying function (i.e.\ $L(ax)/L(x)\rightarrow 1$ as $x\rightarrow\infty$
for every $a>0$) and suppose that the tail $\bar F$ of the service-time distribution satisfies  
$\bar F(x)= L(x) x^{-\alpha}$ for $\alpha \in (1,2)$. Let $c_r$ be the solution of $c_r \bar F(c_r)=1/r$. 
Now, for each $r$, define $\hat Q_r^\gamma(t) = Q_r^\gamma (c_r t)/r, t\geq 0$.  \\

\noindent
{\bf Problem 2:} identify a candidate limit process $Q_*^\gamma(\cdot)$ and show that $\hat Q_r^\gamma (\cdot)\rightarrow Q_*^\gamma(\cdot)$ in distribution as $r\rightarrow\infty$.

\section{Discussion} 
Problem 1 has been solved for PS in  \cite{Gromoll}, \cite{PW} for general (with finite second moment) inter-arrival time distribution and a finite fourth moment assumption on the service-time distribution. The fourth moment assumption was dropped in \cite{LSZ13} using a different method, which required Poisson arrivals. For LCFS, Problem 1 has been solved in \cite{Limic1} for Poisson arrivals and extended to renewal arrivals in \cite{Limic2}. 

The drift and the variance of the RBM for PS and LCFS are identical, which is no coincidence, because it is shown in \cite{FZ11} that the distribution of $Q^\gamma_r (t)$ (for fixed $t$) is independent of $\gamma$, if $Q_r^\gamma(0)=0$. Given this, it is natural to conjecture that the scaling limit in Problem 1 should be independent of $\gamma$. It is remarkable though that this independence of $\gamma$ may break down if the Poisson arrival assumption is relaxed. For the LCFS case, this is illustrated in \cite{Limic2}.

To solve Problem 1, one can think of several strategies. The approach in  \cite{Gromoll}, \cite{PW} is to use a detailed measure-valued description of the queueing dynamics, and hinges on the physical intuition that, in heavy traffic, the queue length becomes a deterministic function of the workload. This property, known as state-space collapse, should still hold for general $\gamma$. The approaches in \cite{LSZ13} and \cite{Limic1} are based on connections with branching processes and excursion theory. To make such an approach work, one would first have to identify an appropriate relation between a Crump-Mode-Jagers (CMJ) branching process (see \cite{LSZ13} and references therein for background on CMJ processes) and the service discipline $\gamma$. A more pragmatic way is to first prove Problem 1 for the special case of phase-type service-time distributions. 

Problem 2 has been haunting me for the past 20 years. I view Problem 2 to be much harder than Problem 1. One reason is that the limiting queue-length process can no longer be a deterministic function of the heavy traffic limit of the workload process (for the latter, see \cite{Whitt2002}). Note that the distinction between finite and infinite variance is not needed for the invariant distribution in the case: when $\rho_r<1$ for all $r$,  $(1-\rho_r) Q_r^\gamma(\infty)$ converges in distribution to an exponential random variable with unit mean as $r\rightarrow\infty$. 

It is possible to determine weak convergence of $\hat Q^\gamma_r(t)$ for fixed $t$ as $r\rightarrow\infty$, assuming $Q_r^\gamma(0)=0$, using the result in \cite{FZ11}. This suggests that marginal distributions of a candidate limit process $Q^\gamma_*(\cdot)$ should be independent of $\gamma$. This will not be the case for the entire process, which has been determined in \cite{Duquesne} for LCFS and in \cite{LSZ13} for PS. In both cases, the limiting queue-length process is expressed in terms of a functional of a sequence of continuous-state branching processes, namely the height process in the LCFS case, and a Lamperti transform in the PS case; see \cite{Duquesne} and \cite{LSZ13} for detailed descriptions. 

For PS, convergence is still an open problem, as tightness of $\{\hat Q_r^{PS} (\cdot)\}_r$ has not been established. In the case where $\gamma$ is completely general, I have no conjecture on what the candidate limit process $Q^\gamma_*(\cdot)$ might be. 

An idea which may lead to partial results of independent interest is to try and establish a relationship between the symmetric queue governed by $\gamma$, and some functional of a CMJ process, which should be dependent on $\gamma$. As mentioned in the discussion of Problem 1, such a connection would be interesting in its own right. Once such a connection has been established, it may be possible to use the technique in \cite{LS14}, leveraging the fact that $\hat Q_r^\gamma (\cdot)$ is a regenerative process for each $r$.

A sufficient condition for tightness  of $\{\hat Q_r^{PS} (\cdot)\}_r$  has been described in \cite{LS15} which is appealing, as it is stated in terms of a tail bound for the distribution of 
the maximum queue length during a busy cycle. Deriving such a tail bound for arbitrary symmetric service disciplines could be an interesting problem in its own right. To develop intuition, one may first approach this by assuming phase-type service times. A related tail bound for bandwidth sharing networks has been derived in \cite{Zhong}.

\bibliographystyle{abbrv}
\bibliography{q18refs}
\end{document}